# The Framework of Mechanics for Dynamic Behaviors of Fractional-Order Stochastic Dynamic Systems[*]

October 8, 2021

**Ruibin Ren[1] and George Xianzhi Yuan[2,3,4,5]**

[1]College of Maths., Southwest Jiaotong Univ., Chengdu 610031 China

[2]Business School, Sun Yat-sen University, Guangzhou 510275 China

[3]Business School, Guangxi Unievrsity, Nanning 530004 China

[4]Business School, Chengdu University, Chengdu 610106, China

[5]BBD Tech. Co.Ltd., #9 Building, Tianfu Avenue 966, Chengdu 610093 China

## Abstract

The goal of this paper is to establish a general framework for dynamic behaviors of coupled fractional-order stochastic dynamic systems of particles by using star-coupled models. In particular, the general mechanics on the dynamic behaviors related to the stochastic resonance (SR) phenomenon of a star-coupled harmonic oscillator subject to multiplicative fluctuation and periodic force in viscous media are established by considering couplings, memory effects, the occurring of synchronization linked to the occurring of SR induced. The multiplicative noise is modeled as dichotomous noise and the memory of viscous media is characterized by fractional power kernel function. By using the Shapiro-Loginov formula and Laplace transform, the analytical expressions for the first moment of the steady state response, the stability and relationship between the system response and the system parameters in the long-time limit in terms of systems' asymptotic stability are also established.

The theoretic and simulation results show the non-monotonic dependence between the response output gain and the input signal frequency, noise parameters provided by fractional-order stochastic dynamics are significant different by comparing those exhibited by traditional integer-order stochastic dynamics, which indicates that the bona fide resonance and the generalized SR phenomena would appear. Furthermore, the fluctuation noise, the number of the particles for the systems, and the fractional order work together, producing more complex dynamic phenomena compared with the traditional integral-order systems.

Finally, all theoretical analyses are supported by the corresponding numerical simulations, and thus it seems that the results established in this paper would provide a possible fundamental mathematical framework for the study of Schumpeter's theory on the economic development under the "innovation and capital paradigm" and related disciplines. In particular, the framework established by this paper allows us at the first time logically concluding that **"in general, the ratio for SMEs grow up successfully is less than one third"**, this is consistent with what the market has been observed commonly, but similar conclusion not available from the existing literature today.

# Keyword

Langevin equation, Star-Coupled structure, Caputo fractional derivative, fractional-Order dynamic system, stochastic resonance (SR), U-Shaping phenomenon, synchronous behavior, MFPT, Schumpeter's theory of economic development.

# 1 Introduction

As the research frontier of the statistical physics and the stochastic dynamical system, the stochastic resonance (SR) driven by fluctuation and periodic signal recently become a popular research direction [2, 3, 6, 5]. The term of SR was proposed by Benzi et al.[7, 8] and Nicolis [9] to explain the climatic mechanics of periodic glaciers in the 1980s. Contrary to the common knowledge that noise is harmful, the SR phenomenon shows that random disturbance (noise) can produce cooperative effect under certain conditions, this means SR can realize the transfer of noise energy to signal energy, and thus strengthen the system output. Since then, more and more scholars have paid attention to the theoretical and experimental researches on SR, which makes it gradually become a hot topic in the field of stochastic dynamics.

Through this paper, the term "SR" means the generalized SR unless a special explanation, which means it is the non-monotonic transformation phenomenon of some functions of system response (such as moment, power spectrum, auto-correlation function, signal to noise ratio, and so on) with some characteristic parameters of the system (such as frequency, excitation amplitude or noise intensity, correlation rate), and more related discussion by Gitterman [10] Yu et al.[76], Gao et al.[77] and references wherein in details.

In the past 30 years, a large part of the research on the SR phenomena have been carried out by different dynamical systems and noise forms, and corresponding physical models have been established respectively. From the perspective of the model, fluctuations enter the model equation in the form of multiplication[13, 14, 6], so the study of the SR phenomenon of the harmonic oscillator system essentially belongs to the study of the dynamic behavior of the resonant subsystem under the multiplicative external noise [11, 12, 13].

Based on the the discussion mentioned above, an intriguing and significant question arises naturally is whether the presence of the potential fluctuations and star-coupling can affect the dynamical properties of overdamped Brownian particles? Regrettably, as far as we know, related research reports are currently rare in existing literature. To elaborate on this question, we will investigate a simple model for coupled Brownian particles moving in a harmonic potential driven by periodic force and the potential fluctuations generated by a dichotomous colored noise trough the fractional stochastic dynamic systems.

The goal of this paper is to discuss the dynamic behaviors for the evolution of coupled fractional dynamic systems for particles by using star-coupled models (for particles) using star-coupled models under the general framework of fractional Langevin equations. We will focus on how the signal "$A_0 \sin \Omega t$" affects the outputs of the harmonic oscillators for the fractional stochastic dynamic systems, through the couplings, memory effect, occurring of synchronization induced the phenomenon of stochastic resonance under the process of systems' asymptomatic stability, by using output gain quantity chosen for describing the dynamic effect in the harmonic medium which is $G :\triangleq \frac{A_{st}}{A_0}$, where $A_{st}$ is the amplitude of the long time stable response, and $A_0$ is the amplitude of the input periodic signal, and $\Omega$ represents the "periodic driving frequency" associated with time $t$ (see more discussion in Section 2 below in details). Moreover, results established in this report might be a possible fundamental framework for the new theory of economics development under the "innovation and capital paradigm" and economics and related disciplines. In particular, the framework established by this paper allows us at the first time logically concluding that **"in general, the ratio for SMEs grow up successfully is less than one third"**, this is consistent with what the market has been observed commonly. Here we like to remark that to the best of authors' acknowledge, it is not available from existing literature for this kind of similar conclusions today; and secondly, as the main goal of this paper is to give a general study for unique features for mechanics of dynamic behaviors for fractional-order dynamic systems by using the star-coupled network structure as the basic model to conduct the study, thus the general conclusion for the evolution of SMEs' growth is given under the fractional-order world for $\alpha \in (0, 0.8]$, not with true data of SMEs from the market, but we will conduct the work in this direction by applying stochastic models with more complex structure considering the true samples of SMEs in the practice in coming research work.

This paper is organized as follows. After the introduction, the system model has been solved analytically in section 2. In section 3, applying analytical results and numerical simulations, we discuss the impacts of system parameters, noise intensity and external driving frequency on the output gain $G$, respectively. In section 4, the general dynamic behaviors for the coupled fractional dynamic systems with focus on the phenomenon of stochastic resonance associated with synchronization are given. Finally, the conclusion with brief discussion is summarized in section 5.

## 2 The Fractional-Order Langevin System for Mean Fields and Particles

### 2.1 The Framework of Fractional Langevin Equations

The classical integer-order Langevin equation inherits and develops Einstein's theoretical research on Brownian motion, and simplifies the fact of a large number of molecular collisions into the random perturbation force of the systems (see [2]-[13] and references wherein), which can be expressed as follows: to begin with, let's consider the classical integer order Langevin equation(see [1])

$$m\ddot{x}(t) = -\gamma \dot{x}(t) + F(x,t) + \eta(t), \tag{1}$$

Here, $m$ represents the mass of the particle, $x(t), \dot{x}(t), \ddot{x}(t)$ are the displacement and the corresponding velocity and acceleration of the time-dependent particle. $m\ddot{x}(t)$ represents the inertia force, $\gamma \dot{x}(t)$ represents the damping force, and $F(x,t)$ represents the deterministic driving force which usually contains the potential field force and the applied driving force.

Here $\eta(t)$ is the internal noise, without loss of generality, we will omit it in this manuscript due to the following two reasons: first, compared with the external driving force $F(x,t)$, $\eta(t)$ is so weak that its effects on the system evolution are negligible; second, we mainly focus on the first moment $\langle x(t) \rangle$, and thus $\eta(t)$ vanishes under the general assumption that $\eta(t)$ has a zero mean.

In the above model (1), the term of damping force is equivalent to the following expression:

$$\gamma \dot{x}(t) = \gamma \int_0^t \delta(t-\tau)\dot{x}(\tau)d\tau. \tag{2}$$

which means that the damping force on the particle is dependent only on its current velocity $\dot{x}(t)$. However, more studies indicate that in heterogeneous media, especially in viscous media and fluids with internal degrees of freedom, most physical and biochemical reaction processes show power-law memory of historical states, which means the closer to the current moment, the stronger its memory and the greater the impact on the system.

Therefore, the following damping kernel function with power-law memory has been introduced:

$$\gamma(t) = \frac{1}{\Gamma(1-\alpha)}|t|^{-\alpha}, (0 < \alpha \leq 1). \tag{3}$$

In Fig.1, with the increase of time $t$, $\gamma(t)$ attenuates as a power function, and the larger $\alpha$ is, the faster the decay rate of $\gamma(t)$ is, which means the poorer memory of the system is reflected. Therefore, the damping force evolution characteristics of particles in viscous media can be well described by the fractional damping kernel function $\gamma(t)$. In addition, on the molecular scale, the Regnold number indicating the ratio of inertia force to viscous force is very small, namely, the effect of inertia on the motor is far less than that of damping, therefore, the inertia term in model (1) can be reasonably ignored, namely, the over-damping case with power-law memory is considered as follows:

$$\frac{1}{\Gamma(1-\alpha)} \int_0^t (t-\tau)^{-\alpha} \dot{x}(\tau) d\tau = F(x,t), (0 < \alpha \geq 1) \tag{4}$$

On the other hand, according to the definition of fractional-order Caputo $\alpha$ derivative (see [55], [56])

$$_0^C D_t^\alpha x(t) = \frac{1}{\Gamma(1-\alpha)} \int_0^t (t-\tau)^{-\alpha} \dot{x}(\tau) d\tau, (0 < \alpha \leq 1), \tag{5}$$

Therefore, Eq.(3) can be written as:

$$_0^C D_t^\alpha x(t) = F(x,t), (0 < \alpha \leq 1) \tag{6}$$

Eq.(4) formally constitutes a fractal-order Langevin equation, which is also called fractal-order Langevin model. This model takes into account the power memory characteristics of non-uniform media. Especially, Eq.(4) degenerates into the classical Langevin equation when the fractional order $\alpha = 1$ (see [55], [56] and references wherein).

In what follows, we recall the definition of the stability for the system of fractional differential equations (see more from [55] and Chapter 5 of [57] in details).

We all know that the stability for the system of (fractional) differential equations is a very important topic, and various notions of stability are commonly discussed (e.g., see Chapter 5 in [57]). In general, the stability issues are usually investigated for first-order equations, i.e., for equations requiring only one initial condition to guarantee the uniqueness of the solution. But for any given $\alpha \in (0,1]$ in the equation (6), if we assume $x(t)$ is a function mapping to $R^n$ for a given integer $n > 0$, and of course $F$ must then be defined on a suitable subset of $R^{n+1}$. When talking about the stability, one should be interested in the behaviour of the solutions of (6) for $t \to \infty$. Therefore we will only consider problems whose solutions $x$ exist on $[0, \infty)$. The first of these assumptions is that $F$ is defined on a set $D := [0, \infty) \times \{w \in R^n : \|w\| \leq W\}$ with some $0 < W < \infty$. The norm in this definition of $D$ may be an arbitrary norm on $R^n$. Our second assumption is that $F$ may be a continuous on its domain of definition and that it satisfies such as a Lipschitz condition as this may assert that the initial value problem consisting of (6) and the initial condition $x(0) = x_0$ has a unique solution on the interval $[0, b)$ with some $b \leq \infty$ if $\|x_0\| \leq W$. Finally we assume that $F(t, 0) = 0$ for all $t \geq 0$, this condition would imply that the function $x(t) = 0$, which is a solution of fractional differential equation (6). Under the motivation for these hypotheses discussed here, we have the concepts for the stability of system governed by (fractional) differential equations as the follows.

**The Definition of Stability for the Fractional Systems:**

(a) The solution $x(t) = 0$ of the fractional differential equation (6), subject to the assumptions mentioned above, is called **stable** if, for any $\epsilon > 0$, there exists some $\delta > 0$ such that the solution of the initial value problem consisting of the differential equation (6) and the initial condition $x(0) = x_0$ satisfies $\|x(t)\| < \epsilon$ for all $t \geq 0$ whenever $\|x_0\| < \delta$.

(b) The solution $x(t) = 0$ of the differential equation (6), subject to the assumption mentioned above, is called **asymptotically stable** if it is stable and there exists some $\gamma > 0$ such that $\lim_{t \to \infty} \|x(t)\| = 0$ whenever $\|x_0\| < \gamma$.

It is clear that the notion of stability is the extension of this idea to unbounded intervals, which means the trivial solution is stable if a small change in the initial value leads to a small change of the solution over the complete positive half-line, obviously this is much stronger than the continuous dependence on the given data. Secondly, asymptotic stability is even stronger since it requires the solution of the perturbed problem not only to remain close to the original solution but actually to converge to the latter.

## 2.2 The Framework of Star-Coupled Fractional-Order Systems

In this paper, we study the collective behavior of $N + 1$ star-coupled fractional harmonic oscillators with fluctuating frequency driven by external periodic force:

$$\begin{cases} {}^C_0D^\alpha_t x_0(t) = (\omega + \xi_t)x_0(t) + \epsilon \sum_{i=1}^{N}(x_i - x_0) + A_0\sin(\Omega t); \\ {}^C_0D^\alpha_t x_k(t) = (\omega + \xi_t)x_k(t) + \epsilon(x_0 - x_k) + A_0\sin(\Omega t), k = 1, 2, ..., N. \end{cases} \quad (7)$$

Here, $x_i, i = 0, 1, 2, \cdots, N$) is the position of the $i$th particle at time $t$. Particularly, $x_0(t)$ and $x_i(t)(i \neq 0)$represent the main particle and the general particles of the star-coupled oscillator, respectively, this is the structure by the right panel of Fig.1 below).

Here, we like to share with readers by a simple reason for the study of star-coupled dynamic system (7) is this case: think of a small or medium-sized enterprise $x_0$ (SME, in short) at its initial stage, it is reasonable to assume that the size $N$ of the SME $x_0$'s partners (i.e., all of related investors or senior management members) should be a smaller number, plus each of them are not (highly) related to each other's in terms of investment or management relationship, this thus could be represented (at least approximately) by the "Star-Network Structure" as shown by the right panel of Figure 1.

The oscillators interact with each other through a linear coupling term $\epsilon(x_0 - x_i)$ or $\epsilon\sum_{i=1}^{n}(x_i - x_0)$, $\epsilon$ is called coupled coefficient, and the whole system is driven by the external periodic force $A_0\sin(\Omega t)$, with $A_0$ and $\Omega$ represent the amplitude and periodic driving frequency, respectively.

For coupled particles with dichotomous noise in the dynamical system, we introduce star coupled particles in a randomly switching potential. Hence, the fluctuation potential can be modeled by the term "$(\omega + \xi_t)x_0(t)$" or "$(\omega + \xi_t)x_i(t)$", here $\omega \geq 0$ is the intrinsic frequency of the system that disturbered by the potential fluctuation $\xi_t$; and also $\xi_t$ is modeled as a dichotomouse Markov process, which also called the random telegraph noise.

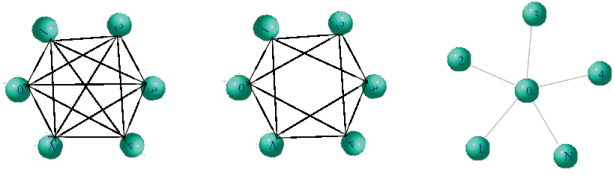

Figure 1: The different structures of network systems

Furthermore, $\xi_t$ is a symmetric dichotomous noise that randomly switches between two values $\pm\sigma(\sigma > 0)$ with the mean value and correlation function:

$$\langle \xi_t \rangle = 0, \tag{8}$$

and

$$\langle \xi_t \xi_\tau \rangle = \sigma^2 \exp(-\lambda |t - \tau|), \tag{9}$$

respectively.

## 2.3 The Framework of Mean Field Modelling for Star-Coupled Fractional-Order Systems

To discuss the collective behavior of the coupled system, we introduce the mean field defined as following formula:

$$S = \frac{\sum_{i=0}^{N} x_i}{N+1} \tag{10}$$

here, $S$ is a profile of the average displacement of all the $N+1$ particles.

By averaging all $N+1$ equations statistically with the mean field $S$, we have

$$D^\alpha S + (\omega + \xi_t)S = A_0 \sin(\Omega t). \tag{11}$$

To get the first moment of $S$ of the steady-state, we firstly average all realizations of the trajectory in Eq.(11), thus we obtain

$$D^\alpha \langle S \rangle + \omega \langle S \rangle + \langle \xi_t S \rangle = A_0 \sin(\Omega t). \tag{12}$$

Since Eq.(12)contains a new correlator: $\langle \xi_t S \rangle$, we established an another equation by multiplying Eq.(11) by $\xi_t$ and averaging the obtained equation, which yields

$$\langle \xi_t D^\alpha S \rangle + \omega \langle \xi_t S \rangle + \sigma^2 \langle S \rangle = 0. \tag{13}$$

Using the well-known Shapiro-Loginov formulas from [54], we have

$$\langle \xi_t D^\alpha S \rangle = e^{-\lambda t} D^\alpha (\xi_t S e^{\lambda t}). \tag{14}$$

Indeed, by the definition of Caputo fractional calculus, we have

$$\langle \xi_t D^\alpha S \rangle = \left\langle \xi_t \int_0^t \frac{(t-u)^\alpha \dot{S}(u)}{\Gamma(1-\alpha)} du \right\rangle = \int_0^t \frac{(t-u)^\alpha}{\Gamma(1-\alpha)} \left\langle \dot{S}(u) \xi_t \right\rangle du \tag{15}$$

using the Shapiro-Loginov formula[54], we have

$$\frac{d \langle \xi_t \dot{x}(u) \rangle}{dt} + \lambda \langle \xi_t \dot{x}(u) \rangle = 0, \tag{16}$$

by solving Eq.(16), it follows that

$$\langle \xi_t \dot{x}(u) \rangle = e^{-\lambda(u-t)} \langle \xi_u \dot{x}(u) \rangle \tag{17}$$

thus the Eq.(14) is established. Inserting Eq.(14) into Eq.(13), we obtain the following expression:

$$e^{\lambda t} D^\alpha (\xi_t S e^{\lambda t}) + \omega \langle \xi_t S \rangle + \sigma^2 \langle S \rangle = 0. \tag{18}$$

Let $y_1 = \langle S \rangle, y_2 = \langle \xi_t S \rangle$, we obtain a closed equations for two variables $y_1, y_2$ according to Eq.(12) and (18):

$$\begin{cases} D^\alpha y_1 + \omega y_1 + y_2 = A_0 \sin(\Omega t), \\ e^{\lambda t} D^\alpha (y_2 e^{\lambda t}) + \omega y_2 + \sigma^2 y_1 = 0. \end{cases} \tag{19}$$

Applying the Laplace transform technic to the closed Eqs.(19), we have:

$$C \begin{pmatrix} Y_1 \\ Y_2 \end{pmatrix} = \begin{pmatrix} s^{\alpha-1} y_1(0) + \frac{A_0 \Omega}{s^2 + \Omega^2} \\ (s+\lambda)^{\alpha-1} y_2(0) \end{pmatrix} \tag{20}$$

Here, $Y_1$ and $Y_2$ are the Laplace transformation of $y_1$ and $y_2$, $y_1(0)$ and $y_2(0)$ are the initial conditions. $C = (c_{ij}), i = 1, 2; j = 1, 2$, and

$$C = \begin{pmatrix} s^\alpha + \omega & 1 \\ \sigma^2 & (s+\lambda)^\alpha + \omega \end{pmatrix} \tag{21}$$

Then the corresponding characteristic equation of (19) is

$$(s^\alpha + \omega)((s+\lambda)^\alpha + \omega) - \sigma^2 = 0. \tag{22}$$

Suppose that $s = \theta i = |\theta| (\cos \pi/2 + i \sin(\pm \pi/2))$ is a root of Eq.(22), where $\theta$ is a real number and when $\theta > 0$, we take $i \sin(\pi/2)$ while when $\theta < 0$, we take $-i \sin(\pi/2)$. Then one has

$$((|\theta| (\cos \pi/2 + i \sin(\pm \pi/2)))^\alpha + \omega)((|\theta| (\cos \pi/2 + i \sin(\pm \pi/2)) + \lambda)^\alpha + \omega) - \sigma^2 = 0. \tag{23}$$

By separating real and imaginary parts, we have

$$\omega^2 - \sigma^2 + \omega R^\alpha \cos(\alpha\Theta) + \omega\theta^\alpha \cos(\alpha\pi/2) + \theta^\alpha R^\alpha \cos(\alpha(\Theta + \pi/2)) = 0 \quad (24)$$

and

$$\omega R^\alpha \sin(\alpha\Theta) + \theta^\alpha(\omega \sin(\alpha\pi/2) + R^\alpha \sin(\alpha(\Theta + \pi/2))) = 0. \quad (25)$$

where $R = \sqrt{\lambda^2 + \omega^2}, \Theta = \arctan(\frac{\theta}{\lambda})$. By two equations above, we have

$$(\omega^2 - \sigma^2 + \omega R^\alpha \cos(\alpha\Theta) + \omega\theta^\alpha \cos(\alpha\pi/2) + \theta^\alpha R^\alpha \cos(\alpha(\Theta + \pi/2)))^2 + \\ (\omega R^\alpha \sin(\alpha\Theta) + \theta^\alpha(\omega \sin(\alpha\pi/2) + R^\alpha \sin(\alpha(\Theta + \pi/2))))^2 = 0. \quad (26)$$

Obviously, Eq.(26) has no real solutions, which means that Eq.(22) has no purely imaginary roots under the system assumption. Thus if we choose

$$\sigma^2 < \omega^2 + \omega\lambda^\alpha \quad (27)$$

then two eigenvalues of coefficient matrix $C$ of (19) are negative.

Now by Theorem 7.20 for the stability of linear fractional differential system [55](see also [58, 59]), all the root $s$ of the characteristic equation satisfy $|arg(c)| > \alpha\pi/2$, we obtain the necessary and sufficient condition of the stability criteria of Eq.(19) given by Eq.(27). Especially, when $\alpha = 1$, Eq.(27) will degraded into the stability criterial of the integer order system.

Solving the Eqs.(20), we have the following

$$Y_1(s) = H_1(s)y_1(0) - H_2(s)y_2(0) + H_0 \frac{A_0 \Omega}{s^2 + \Omega^2}, \quad (28)$$

where

$$\begin{aligned} H_0 &= \frac{c_{22}}{c_{11}c_{22} - c_{12}c_{21}}, \\ H_1 &= \frac{-c_{21}}{c_{11}c_{22} - c_{12}c_{21}} s^{\alpha-1}, \\ H_2 &= \frac{-c_{12}}{c_{11}c_{22} - c_{12}c_{21}} (s+\lambda)^{\alpha-1}. \end{aligned} \quad (29)$$

Solving the equation (28) and applying the inverse Laplace transform to Eq.(20),

$$\langle S \rangle = y_1 = A_0 \int_0^t h_0(t-\tau) \sin(\Omega\tau) d\tau. \quad (30)$$

Here, $h_0$ is the inverse Laplace transform of $H_0$.

Therefore, the first moment $\langle S \rangle$ of the stable-state mean field system can be considered as the output of the linear time-invariant system (11) with the transform function $h_0(t)$. On the other hand, according to the response theory of linear time-invariant system, the output can be written as:

$$\langle S \rangle = A_1 \sin(\Omega t + \phi). \quad (31)$$

Where $j$ is the imaginary unit which means $j^2 = -1$, $A_1 = A_0|H_0(j\Omega)|$ and $\phi$ are the output amplitude and the phase angle of, respectively. Where,

$$|H_0(j\Omega)| = \sqrt{\frac{f_1^2 + f_2^2}{f_3^2 + f_4^2}}, \phi = \arctan(\frac{f_2 f_3 - f_1 f_4}{f_1 f_3 + f_2 f_4}). \tag{32}$$

Finally, we obtain the output gain (OAG) for the mean field of star-coupled dynamic system:

$$G = \frac{A_1}{A_0} = |H_0(j\Omega)| = \sqrt{\frac{f_1^2 + f_2^2}{f_3^2 + f_4^2}}. \tag{33}$$

Where the explicit expression of $f_i, (i = 1, 2, 3, 4)$ is given below:

$$\begin{aligned}
f_1 &= \omega + r^\alpha \cos(\alpha\theta), \\
f_2 &= r^\alpha \sin(\alpha\theta'), \\
f_3 &= r^\alpha \Omega^\alpha \cos(\alpha\theta + \frac{\alpha\pi}{2}) + \omega\Omega^\alpha \cos(\frac{\alpha\pi}{2}) + \omega f_1 - \sigma^2, \\
f_4 &= r^\alpha \Omega^\alpha \sin(\alpha\theta + \frac{\alpha\pi}{2}) + \omega\Omega^\alpha \sin(\frac{\alpha\pi}{2}) + \omega f_2, \\
r &:= \sqrt{\Omega^2 + \lambda^2}, \\
\theta &:= \arctan(\frac{\Omega}{\lambda}).
\end{aligned} \tag{34}$$

It should be noted that the results in this section recovery the previously published results for a single oscillator without any memory effect when both $\epsilon = 0$ and $\alpha = 1$.

## 2.4 Modeling Star-Coupled Fractional-Order Systems

To analyze the SR behavior of the coupling system, firstly, we need to know whether or not the average behaviors of the $N + 1$ particles are synchronous. We know from the previous section that the stable response of the mean field is a sinusoidal wave, thus the synchronization of the single particle are further considered, including both the main particle and the general particles in the coupled system.

### 2.4.1 Modelling Synchronization between Mean Fields and Main Particle

Firstly, we consider the main particle which is governed by the following equation:

$$D^\alpha x_0 + (\omega + \xi_t)x_0 = \epsilon \sum_{i=1}^{N}(x_i - x_0) + A_0 \sin(\Omega t). \tag{35}$$

Denoting the deviation of the main particle's displacement from the mean field as: $\Delta_0 = x_0 - S$, using Eq.(35) minus Eq.(11) yields

$$D^\alpha \Delta_0 + (\omega + \xi_t)\Delta_0 = \epsilon \sum_{i=1}^{N}(x_i - x_0). \qquad (36)$$

Notes that,

$$\Delta_0 = x_0 - \frac{1}{N+1}\sum_{i=0}^{N} x_i = \frac{1}{N+1}\sum_{i=1}^{N}(x_0 - x_i), \qquad (37)$$

Then we have

$$\sum_{i=1}^{N}(x_i - x_0) = -\Delta_0(N+1). \qquad (38)$$

Inserting Eq.(38) into Eq.(36) and simplifying the observed equation, we have

$$D^\alpha \Delta_0 + [\omega + \xi_t + \epsilon(N+1)]\Delta_0 = 0. \qquad (39)$$

With the same method in section 2.3, we obtain the condition of the stability criteria of (36):

$$\sigma^2 < [\omega + \epsilon(N+1)]^2 + \lambda^\alpha[\omega + \epsilon(N+1)], \qquad (40)$$

Since Eq.(39) is a homogeneous equation, it must has a zero solution which is globally and asymptotically stable when Eq.(39) satisfied the condition of the stability criteria Eq.(40). So when $t \to \infty$, we have

$$\langle \Delta_0 \rangle = 0. \qquad (41)$$

This indicate that the main particle will synchronize with the mean field of the system in expectation when the evolution time is long enough. In this paper, Eq.(40) is called the synchronous criteria of the main particle. This is one of the main results of this study.

### 2.4.2 Modelling Synchronization Among Mean Fields and General Particles

In the next section, we further consider the general particle of system Eq. (1), say the $n$th, which is governed by the following equation:

$$D^\alpha x_n + (\omega + \xi_t)x_n = \epsilon(x_0 - x_n) + A_0 \sin(\Omega t), n = 1, 2, ..., N, \qquad (42)$$

Let $\Delta_n = x_n - S$, subtracting Eq.(42) from Eq.(11) we observed,

$$D^\alpha \Delta_n + (\omega + \xi_t)\Delta_n = \epsilon(x_0 - x_n). \qquad (43)$$

When the synchronous criteria of the main particle Eq.(40) is satisfied, $x_0 = S$ a.e., as $t$ goes to infinity. Therefore, Eq.(43) can be rewritten as:

$$D^\alpha \Delta_n + (\omega + \xi_t + \epsilon)\Delta_n = 0. \qquad (44)$$

Comparing with Eq.(39), we obtain the first moment of the general particle's deviation in the long time limit $t \to \infty$:

$$\langle \Delta_n \rangle = 0, \tag{45}$$

and the corresponding synchronous condition of the general particle:

$$\sigma^2 < (\omega + \epsilon)^2 + \lambda^\alpha (\omega + \epsilon). \tag{46}$$

Since when Eq.(46) holds, Eq.(40) must be hold. So we call Eq.(46) as the global synchronization condition, and let

$$\sigma_s^2 = (\omega + \epsilon)^2 + \lambda^\alpha (\omega + \epsilon). \tag{47}$$

When the system satisfied the global synchronous condition, all the particles, including the main particle and the general particles, will be synchronous with the main field of the system $S$ as $t$ becomes larger than a certain threshold.

Furthermore, the results Eq.(33),Eq.(41) and ((45)) tells us that the mean field is equal to the average of any single particle's position in some certain threshold. This allows us to study dynamic behaviors of systems among single particles through the the mean field, related to synchronous, mean first passing time (MFPT), occurring of stochastic resonance, and significant difference on the output gain $G$ for fractional dynamic systems by comparing with traditional dynamic systems with order of the derivative $\alpha$ being or closing the integer 1 in sections 3 and 4 below.

## 3 The Simulations of Behaviors for Star-Coupled Fractional-Order Systems

The goal of this section is to discuss some significant different features for dynamic synchronous behaviors of star-coupled fractional-order systems by comparing with traditional integer order system. It is well known that Stochastic Resonance (SR) would be induced (derived) by the synchronization effect between a noise-induced switch and weak signal modulation in the time-modulated bistable model [30]. Indeed, in the parallel to the big boom of chaos synchronization ([31]-[32]), frequency and phase synchronization in stochastic systems have repeatedly received considerable interest by [33] and several decades after the seminal work of Stratonovich in [35], the revival of effective phase synchronization explains from the discovery that the phenomenon of stochastic resonance (SR) can be reinterpreted in terms of a noise-induced phase synchronization (NIPS) (see [30]-[36], and also the work on the first experimental and numerical observation given by [37]-[38]). At the same time, the relation between aperiodic SR and NIPS with stochastic signals [41] is also established in [39] and [40], which has a direct application in the context of behavioural biology by [44]. Moreover, the phenomenon of noise-induced phase synchronization is a much more stringent effect than SR as shown by Freund et al. in [33], and Al-Khedhairi also established results for the study on the stability and synchronization of chaotic finance models (see [50] and references wherein).

On the other hand, by following the concept "mean first passage time" (MFPT) for dynamic systems, Kang et al. recently established the criteria for the occurring of SR in a transcriptional regulatory (integer-order) systems with non-Gaussian noise in [51] for integer-order. By the fact that the concept of MFPT could be used to describe how fast to reach the (asymptotic) stability for the fractional-order system during its evolution process when the time $t \to \infty$, and star-coupled fractional-order system exhibits "an inverse SR-related effect" with size $N$ of particles, which normally does not occur by systems with integer-order, this kind of significant different behaviors for fractional-order would indicate that with the concept of (asymptotic) stability for fractional stochastic dynamical system, the MFPT for fractional-order system may play an important role for the explanation of mechanics for the evolution and synchronization behaviors of star-coupled fractional-order dynamics in the practice for what we will discuss here.

In what follows, the numerical simulations for fractional-order systems are presented to check whether theoretical results (Eq.(33), Eq.(41) and Eq.(45)) are in accordance with numerical simulations. Then the functional relationship of the MFPT on the size $N$ with different fractional order coupled particles will be discussed subsequently.

For the numerical simulation, we applied the fractional predictor-corrector method, and the simulation parameters are: simulation duration $t = 15$, time step $dt = 0.01$, simulation time $T = 3000$; and also we assume that the initial positions of all the particles in the system obey the normal distribution with mean zero and standard deviation one.

## 3.1 The Synchronous Behaviors of Fractional-Order Dynamic Systems

To analyze the dynamical behavior of fractional-order coupling system (7), firstly, we need to check whether or not the average behaviors of size $N+1$ of particles are synchronous as we deduced in Eq(41) and Eq.(45).

From Fig.2(a) to Fig.2(d), all the red solid lines depict the theoretical trajectory of the first moment of the mean field, which is a sine wave that expressed by Eq.(31). we see all the realizations of the deviation are zero in the long-time regime except for a set of measure zero, which indicates that every single particle of Eq.(1) is almost synchronized with the mean field S after some certain time. Furthermore, from Eq.(37), Eq.(41) and Eq.(35), we can infer that particles with different initial conditions will move together with the mean field in the stationary regime.

In Fig. 2(a), the solid line depicts the theoretical trajectory of the mean field while the dished line with circles depicts the corresponding numerical trajectory of the mean field. The figure shows that the theoretical results and the numerical results of the mean filed are statistically uniformity which indicates that the effectiveness of our analysis.

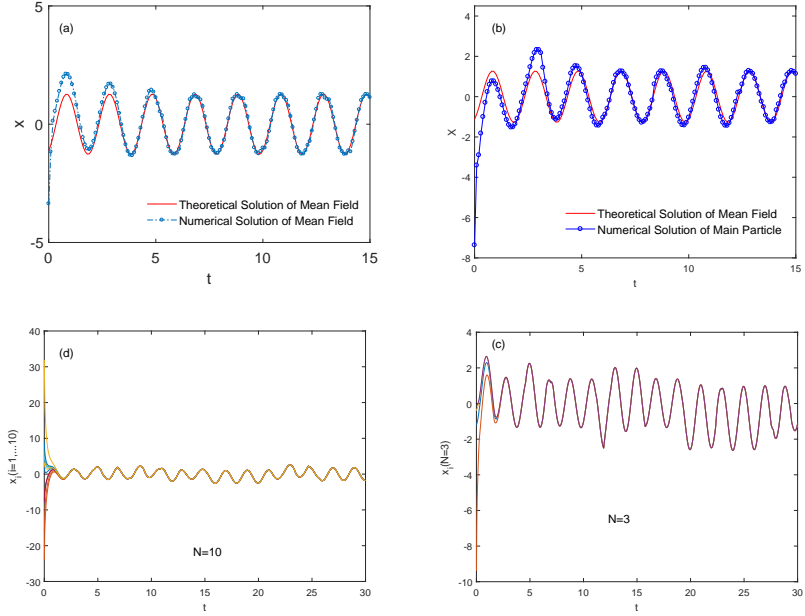

Figure 2: The trajectory of the mean field, the main particle and the general particles. Here, the parameters are set as $\alpha = 0.9, \omega = 1, \epsilon = 1, \lambda = 1, \Omega = \pi, \sigma^2 = 0.5 * \omega^2 + \omega * \lambda^\alpha$. (a):the numerical and theoretical results of the mean field, $N = 10$;(b): the numerical simulation of the main particle vs the theoretical result of the mean field, $N = 10$; (C): the numerical simulations of the general particles vs the theoretical result of the mean field for $N = 3$; and (d): the numerical simulations of the general particles vs the theoretical result of the mean field, $N = 10$.

From Fig. 2(a) to Fig.2(c), different colored lines indicate different displacements of the main and the general particles, the simulation results show that the main particle and every general particles synchronizes with the mean field after a certain time. Here, we mention that the coupling only affects the relaxation time for synchronization(see Fig.2(c) and Fig.2(d)). The essential reason for this phenomenon is that the uniform influence of the randomly switching potential, which can be inferred from inequality Eq.(13) by using a zero-coupling intensity.

As a result, in the stationary regime, we can analyze the dynamical behaviors of all of the particles just by analyzing the mean field for fractional-order systems like the traditional integer-order systems (see [64] and reference where for more in discussion).

## 3.2 The U-Shaping Phenomenon by MFPTs for Star-Coupled Fractional-Order Systems

The mean first passage time (MFPT) which represents an mean initial time that the unstable state transmit to the stable state, is of great importance and may give interesting results concerning the number of coupled particles in the fractional system for the synchronization of particles in the time dimension. Thus the significant role played by first encounters have made calculating the mean first passage time (MFPT) of great importance in numerous real world dynamic systems that can be modeled as random walks (see [60], [61], [62] and references wherein). The first passage time (FPT) is defined as the time required for a random walker to first reach a predefined target state or set of states, however, for any random walk, the exact solution to finding FPT is almost impossible especially in unbounded domain ([63]), and the MFPTs to a target node is the average time expected for a random walker to first find that target on a given network.In general, it is different to have the analytic solutions of MFPTs for a given dynamic system. In this part, based on the fractional-order system (7), we will discuss the influence of coupling on the fractional-order system with simulations.

It is known that the coupling coefficient and the number of the particles are two main factors that influence the coupled system, though the stability conditions of the mean field are independent of the coupling coefficient and the number of the particles (see [64] and related references wherein), but these two factors truly influence synchronous behavior of the system in terms of MFPTs with size $N$ of particles under different kinds of noises. But what we like to address here is that based on the numerical simulations for the coupled fractional-order system, it exhibits significant different behaviors by comparing with integer-order systems, i.e., the numerical results show that the MFPTs are not a simple monotone function of size $N$ any more under the fractional-order system, this is true different from what we observed for integer-order systems (see [64] for more in details), with one important phenomenon exhibited as "an inverse SR-related effect", resulted from the memory effect in supporting "noise enhanced stability" (NES) taking place, which would help us to identify the best range of sizes for systems to reach their synchronization in explaining the mechanics of coupled fractional-order systems' evolution during the process when time $t \to \infty$.

In this section, by numerical simulation, we focused on the MFPTs vs the number of the particles $N$ while varying the value of the fractional order $\alpha$ of systems as shown by (7). Meanwhile, the values of simulations are the same as those used in section 3.1 above.

In Fig.3, the numerical results for MFPTs vs size $N$ shows that the smaller the fractional-order $\alpha$ (here for $\alpha = 0.4$ comparing with $\alpha = 0.9$), the more clearly the so-called "U-Shaping phenomenon", which has obvious non-monotonic relationship with the parameter $N$, and a minimum (valley) can be observed for $N$ somewhat smaller than 12, which turns out to be an "inverse SR-related effect" indicating that the "noise enhanced stability" (NES) phenomenon take-ing place! In most of the cases (e.g., $N > 15$) we have obtained, no minimum can be observed which indicates that the NES effect disappears with the given parameters. Besides, we noticed that the MFPTs in the case of $\alpha = 0.4$ is small than the value in the case of $\alpha = 0.9$, which indicates that the memory effect will enhance the system to arrive the stable state. These finding may be used in effect for the purpose of controlling the MFPTs with different fractional order.

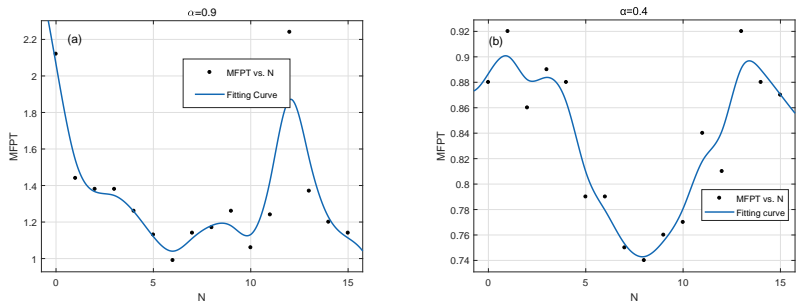

Figure 3: The performance of MFPTs vs particles size $N$. The parameters are: $\omega = 1, \epsilon = 1, \lambda = 1, \Omega = \pi, \sigma^2 = 0.5 * \omega^2 + \omega * \lambda^\alpha$. The left panel: $\alpha = 0.9$; right panel: $\alpha = 0.4$ for the better occurring of "U-shaping phenomenon" by MFPTs for $N$ in the range from 5 to 10 by comparing with one with $\alpha = 0.9$" (left panel).

Here, we like briefly address that the associated numerical results obtained above on the behaviors of MFPTs for star-coupled fractional-order dynamic systems (as shown by Figure 3) are significant different from ones shown by integer-order systems! This new feature for fractional-order dynamic systems of order $\alpha$ (in particular, when less than 0.5) shown as the typical "U-Shaping phenomenon"(equivalently, explained as an "inverse SR-related effect")[1] would provide a possible fundamental mathematical framework for the explanation of new theory for development economy (of enterprises growths with technology innovation) and related disciplines in supporting the economic theory of development due to technology innovations originally introduced and established by the famous economist Schumpeter [73] last century in 1910's.

Indeed, by taking the star-coupled dynamic system in (7) with number size $N$ of particles associated with the main particle $x_0$ as the enterprise $x_0$'s growth from the beginning as one typical small (or medium) enterprise (denoted by SME, in short) with the size $N$ of associated parties (which could be regarded as the size of the so-called "related parties", or say, partners to the enterprise $x_0$ discussed by Yuan et al.[74] and related references) starting with a small number $N$ becoming bigger and bigger during the time period by $t \to \infty$, the process for the growth of the enterprise $x_0$ with the size $N$ of partners should go through at least following three phases.

---

[1] which indicates the capability to "transform disordered information into ordered one" for complex stochastic systems.

**The Phase I: To get the start.** When enterprise $x_0$ is beginning in the SMEs level, the size $N$ normally is very small (zero to around 5), during this time period, the enterprise $x_0$ is facing challenging to get survival by expecting the supporting of synchronous behaviors from size $N$ partners to grow up as discussed by Phase II below !

**The Phase II: To grow up quickly.** By following Phase I, the SME $x_0$ wishes to grow up with a bigger number size $N$ of partners with synchronous action, which means that $x_0$ is trying to reach the stable stage. The numerical results from the right panel of Fig 3 actually show that by choosing the fractional-order $\alpha = 0.4$, the company $x_0$ takes less time to reach synchronization for all partners with the size $N$ being the range from 5 to 10, for this feature by MFPT, which means that $x_0$ with size $N$ is best in the stationary regime by the definition of the stability for dynamic systems, i.e., the system $x_0$ with size $N$ is stable.

**The Phase III: To get the mature by growing strongly.** Once the enterprise $x_0$ with the size $N$ of partners through the Phase II in the stationary regime, to get mature by growing up strongly subject to the emergence of "U-Shaping Phenomenon" (induced) by stochastic resonance under the (asymptotic) stability during a reasonable length of time (thinking of the time $t \to \infty$ ideally).

Putting all three Phases together, it is clear that the enterprise $x_0$ truly grows up unless it goes through Phases I, II and III above, which could be equivalently explained by the way that the company $x_0$ reaches to the stable state from beginning under the supporting of synchronous behaviors by all of its partners with size $N$, and then gets mature through the growth strongly under the occurring of the so-called "U-Shaped Phenomenon" (induced by the emergence of stochastic resonance), which indicates the company $x_0$ in its best state with making of optimal equity and capital structure based on the dynamical evolution under the fractional-order Langevin system framework for SMEs' development (in particular, for the growth/development under the "innovation and capital paradigm" as discussed by Yuan et al.[74], Li et al.[75], Yu et al.[76], Gao et al.[77] for the case of integer order being one).

In addition, the numerical results given by Figures 4 and 5 below in section 4, using star-coupled fractional-order systems as tools to model the mechanics for the growth of the enterprises (denoted by) $x_0$ with a small number $N$ of partners at beginning, during the time for $t$ goes to long enough (thinking of $t \to \infty$), the ratio for the SME $x_0$ growing up successfully (subject to occur SR phenomenon in the stationary regime) is less than 30% by using the fractional-order $\alpha \in (0, 0.8]$ (in particular, when $\alpha$ takes the value around 0.4 is the best). This actually provides a reasonable explanation for the general common observation on the growth of small and medium-sized enterprises (SMEs) with a logical conclusion: "**in general, the ratio of SMEs could grow up successfully is normally not more than one third (actually less than** $30\%$**)**", to the best of our knowledge, this is first time to have the conclusion by using fractional-order stochastic dynamic system as a tool, which is consist with what the market has been observed commonly, not available from any existing literature by so far.

# 4 The Behaviors of Stochastic Resonance for Star-Coupled Fractional-Order Systems

In this section, based on the star-coupled fractional-order system (7) as the tool, we will discuss significant features for the criteria for the classification of stability and the occurring of stochastic resonance phenomena induced by the collective behaviors of the coupled fractional-order systems with numerical simulations respect to output gain $G$ in 3-dimensional plot as functions of the noise intensity $\sigma^2$, the frequency of the system, and external signal "$A_0 \sin \Omega t$" for the conventional stochastic resonance (CSR) and the conventional stochastic resonance (CSR).

Based on the previous theoretical results and numerical analyses, we will firstly discuss the stability criteria, then the SR behaviors will be discussed in detail.

## 4.1 The Stochastic Resonance for Star-Coupled Fractional-Order Systems

As discussed in section 2, the general criteria for the stability of star-coupled fractional-order dynamic systems are that all the roots of $|D(\alpha)| = 0$ have negative real part. These requires the system parameters to satisfy stability conditions Eq.(27).

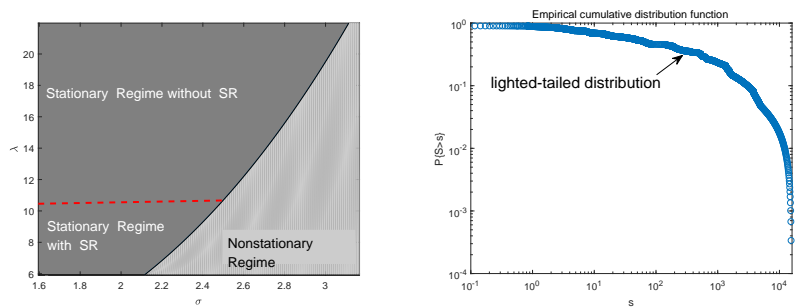

Figure 4: The classification of stationary regime with SR (left panel) by choosing fractional-order $\alpha = 0.8$; and the empirical cumulative distribution function satisfied with time $T = 10$ with time step using $dt = 0.01$ (right panel)

Here, the stability condition involves the calculation of the stability criteria as those established for fractional-order systems in [55] (see also [59], [58] and references wherein), the analytical expressions will be rather cumbersome, therefore, we intend to provide a phase diagram of the $\lambda - \sigma$ plane to analyze the emergence of the stability and the SR.

It is easy to understand that if the system is not stable, it is hard to have a stable state to maintain the status for the occurring of SR, and integer-order stochastic dynamic system one with simple white noise never has SR phenomenon (as the criteria for the SR is not satisfied, see more discussion in [72] and related reference in details). In what following, we discuss systems for which all groups of parameters satisfy the stable condition.

Based on the star-coupled fractional-order dynamic system given by (1) (or (7)) with fractional-order $\alpha = 0.8$, numerical results from the left panel of Fig. 4 shows that the ratio for stable systems with occurring of SR phenomena is around 20%. The parameter region is plotted in Fig.4 with $\omega = 1, \Omega = 1$. The black solid line is the critical line between the nonstationary regime and the stationary regime, while the red dashed line is the critical line for the emergence of SR in the mean field in the stationary regime. Fig.4 shows that SR only occurs in the weak-noise region (the dark grey region with SR). In the strong-noise region (the light grey region with nonstationary response), random forces $\xi_i, i = 1, 2, \cdots, N$ play a dominant role in the system evolution and thus introduce strong randomicity to the system output response. For a fixed time, the $\alpha$ region for the emergence of SR of the domain field by the system (1)(or system (7)) is given in Fig.4 (left panel); and right panel from Fig.4 shows a power-law-like distribution with a light tail which indicates the effectiveness of the fractional system. By the performance of SRs based on fractional-order $\alpha \in [0, 1]$, the results provided by Figure 5 below show us that "smaller the $\alpha$ is, the more likely occurring of SR for the system is"!

Indeed, by using the criteria for the emergence of stochastic resonance (SR) below for the system (1) (or (7)) by the inequality condition (48) below, $0 < r^\alpha \Omega^\alpha \cos(\alpha\theta + \frac{\alpha\pi}{2}) + \omega\Omega^\alpha \cos(\frac{\alpha\pi}{2}) + \omega f_1 < \omega^2 + \omega\lambda^\alpha$, and by choosing the value of fractional-order $\alpha$ in the range $[0, 1]$, the similar numerical results as shown by the performance of SR occurring in Fig. 5 below to confirm that in general the average ratio of systems to occur SRs in the stationary regime is around 28% for the fractional-order $\alpha \in (0, 0.8)$ in general.

Secondly, the comprehensive numerical results by Fig 6 show that based on $G - \Omega - \sigma$ plots for $\alpha = 1, 0.8$ and $0.4$, it shows that the smaller of fractional-order $\alpha$, the stronger the output gain $G$, which is more significant. Here we also like to mention that the output gain $G$ established in this paper would be a useful tool to conduct risk assessment for enterprises' performance by modelling the SMEs and their growths. This unique special feature on the output gain $G$ based on fractional-order less than one would provide another way for us to study the mechanics for the growths of enterprise $x_0$ with size $N$ of partners under the framework of fractional-order Langevin equations by comparing with those did by [74] and [75] for the case of integer-order being one!

## 4.2 The Occurring Criteria of SR Phenomena of Star-Coupled Fractional-Order Systems

This section focus on the synergistic effect of the fractional order $\alpha$ and the noise intensity $\sigma$ and noise correlation rate $\lambda$ on SR behavior, including the bona fide

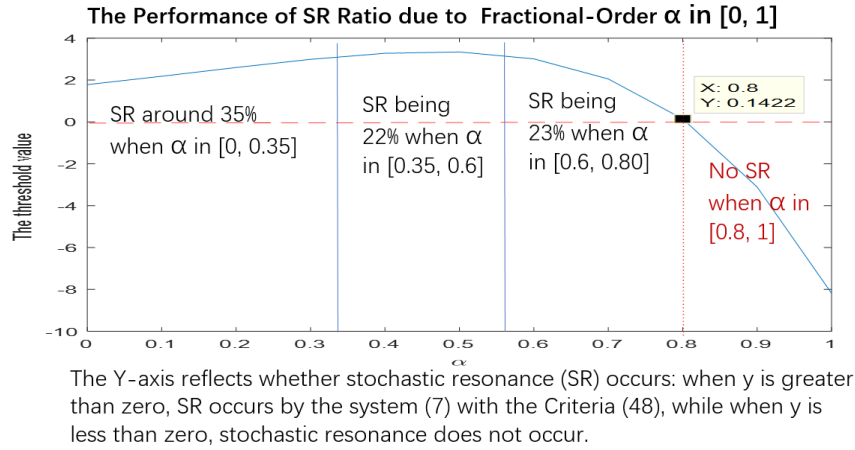

Figure 5: The Ratio for occurring of SRs based on fractional-order taking values in $[0, 1]$

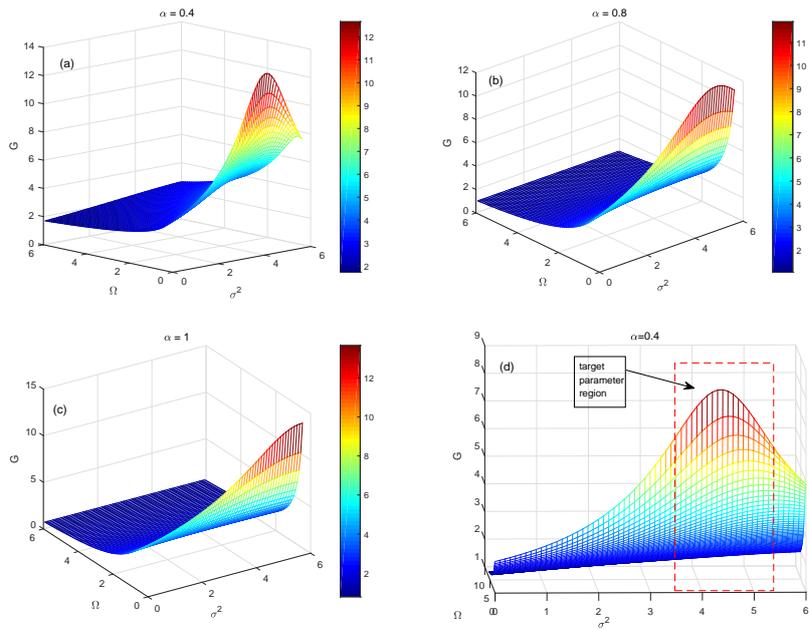

Figure 6: The 3-dimensional plot of the output gain $G$ as function of the noise intensity $\sigma^2$ and the frequency of the periodic sinusoidal excitation $\Omega$; other parameters are set as $\omega = 1, \lambda = 5, A_0 = 1$.

resonance (BSR) (see [71]) by adjusting system parameters with a fixed noise,

or the conventional stochastic resonance (CSR) presented by adjusting the noise intensity with fixed system parameters.

Since the output amplitude $G$ achieves its extreme value, based on Eq.(30) above, we have that when $\sigma^2 := r^\alpha \Omega^\alpha \cos(\alpha\theta + \frac{\alpha\pi}{2}) + \omega\Omega^\alpha \cos(\frac{\alpha\pi}{2}) + \omega f_1 + \omega^2$, it follows that the criterion for the occurring of stochastic resonance (SR) as follows:

$$0 < r^\alpha \Omega^\alpha \cos(\alpha\theta + \frac{\alpha\pi}{2}) + \omega\Omega^\alpha \cos(\frac{\alpha\pi}{2}) + \omega f_1 < \omega^2 + \omega\lambda^\alpha \qquad (48)$$

We like to point out that although the criteria (48) above have nothing to do with the coupling, the coupling intensity or the number of particles determines how fast the synchronization state is reached.

### 4.2.1 The Bona fide Stochastic Resonance (BSR)

In Fig.(7), with the system frequency $\Omega$ increases, some non-monotonic behaviors of $A_1$ have been presented which indicates that the BSR phenomena occurred. Increasing the noise correlation coefficient $\lambda$ while remain the other parameters unchanged we will find the one-peak and one-valley SR which were rarely reported in the single particle system. Moreover, when $\lambda = 2$, the BSR disappeared, the effect comes from the fact that since $\lambda$ increase, the correlation time decrease, then the oscillators will constantly jump from one state to the other, the time to build up a proper response to an external field might be insufficient, therefore, the SR disappears. In general, the case of a coupled oscillator with fluctuation frequency can contribute to the existence of classical BSR, but it can also contribute to the one-peak and one-valley BSR(when the system parameters and the noise parameters are properly matched).

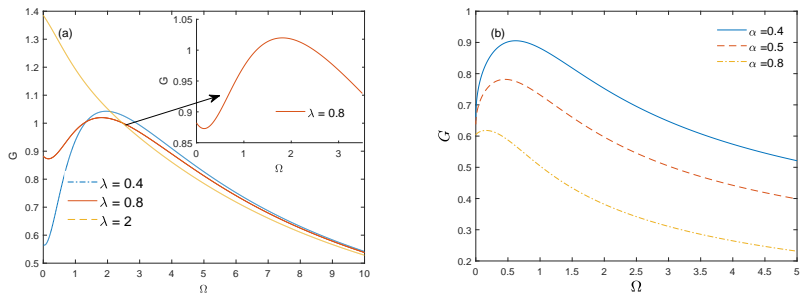

Figure 7: The output gain $G$ as function of $\Omega$ under different value of $\lambda$ and $\alpha$, respectively. The common parameters are set as: $\omega = 1, \sigma = 1$. Other parameters are: (a) $\alpha = 0.9$; and (b) $\lambda = 0.8$.

### 4.2.2 The Conventional Stochastic Resonance (CSR)

We provide the $A_1 - \sigma$ in Fig.(8) to analyze the influence of the noise intensity on the CSR under different $\omega$ and $\Omega$. Based on the conclusions obtained from Eq.(48), Fig.8(a) shows that as $\omega$ increases, the resonance curves switch to the right hand, the peak value decreases, which indicates that increasing $\omega$ can weaken the resonance intensity. From Fig.8(b), the resonance curve decreases almost linearly with the increase of the periodic driving frequency $\Omega$, namely, no CSR occurs. When $\Omega = 0.1\pi, 0.5\pi$, the curves exhibit a resonance peak. Moreover, the resonance peak decreases as the increase of $\Omega$, which means increasing $\Omega$ will weaken the resonance intensity under some conditions. On the other hand, the above results are consistent with the result shown in Fig.(7), that is, changing $\sigma$ can control the emergence of CSR.

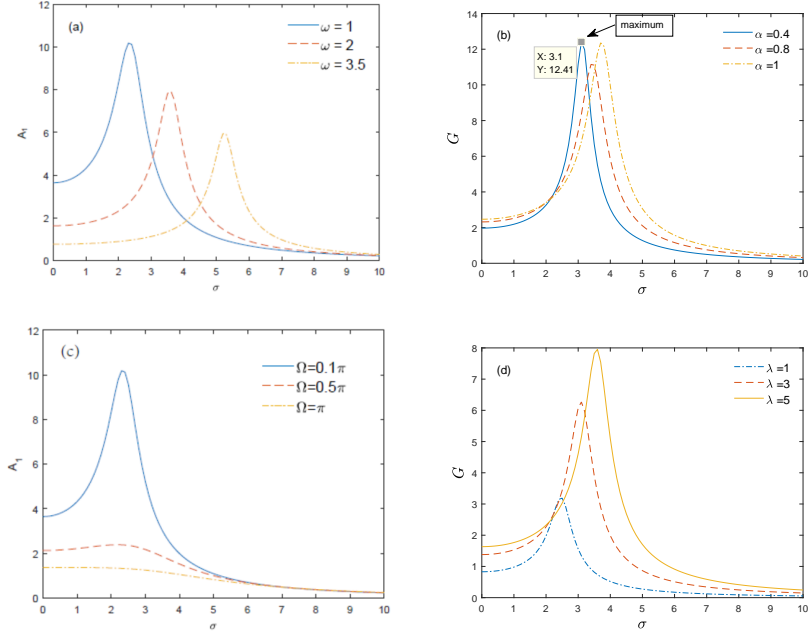

Figure 8: The output gain $G$ as function of $\sigma$ under different value of $\omega$, $\alpha$, $\Omega$ and $\lambda$ ,respectively. Other parameters are: (a) $\alpha = 0.9, \Omega = 0.1\pi, \lambda = 5$; (b) $\omega = 2, \Omega = 0.1*\pi, \lambda = 5$; (c) $\omega = 1, \alpha = 0.9, \lambda = 5$; and (d) $\omega = 2, \alpha = 0.9, \Omega = 0.1*\pi$
.

## 5 The Conclusion with Brief Discussion

To summarize, we investigated, in the long-time regime by considering the asymptotic stability, the collective behaviors including the synchronous and the stochastic resonance phenomenon in a coupled stochastic oscillator with frequency fluctuation driving by periodic signal. The coupling form is the star-coupled oscillator, and the frequency fluctuation means that the potential energy has been perturbed and it is modeled as a dichotomous noise in this paper, and then a closed system will be deduced. By adjusting some parameters of the model, we can effectively control the SR phenomenon of fraction-order coupled systems within a certain range, which provides another new perspective for the application of stochastic resonance theory in practical physics and engineering.

With combination of sections 3 and 4, based on three Phases discussed in section 3 and numerical results for the classification of stationary regimes with occurring of SR given by the left panel in Fig.4, and results provide by Fig.5 for all $\alpha \in [0,1]$, they tell us that by using the fractional-order system (7) as a tool with $\alpha in(0, 0.8])$, the numerical simulation results obtained in sections 3 and 4 provide a fundamental mathematical framework to support us at the first time logically concluding that "**in general, the ratio for SMEs grow up successfully is less than one third (actually less than** $30\%$)", which is consistent with what the market has been observed.

Moreover, we like to briefly address that as shown by the performance of MFPTs with size $N$ given by the right panel of Fig.3, it indicates that the fractional-order (with $\alpha < 0.5$) seems provide a better way to identify the best range of size $N$ for partners (around 5 to 10) for company $x_0$ to have the so-called "U-Shaping Phenomenon" to occur during its growth, which means that the company is best for its making of optimal equity and capital structure based on the dynamical evolution under the fractional-order Langevin system framework for SMEs. Moreover, comparing with results given by Yuan et al.[74], Li et al.[75], Yu et al.[76] and Gao et al.[77] under the framework of stochastic dynamic system with the integer-order one, the theoretic and simulation results obtained in this paper show a number of unique features for fractional-order dynamic systems, which truly significant different from thoese by using $\alpha = 1$! We do expect that the framework for dynamic behaviors of star-coupled fractional dynamic systems established in this paper would help us to provide a new tool for the study of mechanisms in supporting new theory for economic development for SMEs' growth and related topics (in particular, for the development of SMEs under "innovation and capital paradigm" originally introduced by the famous economist Schumpeter [73] last century in 1910s.

Finally, we like to point that as the goal of this paper is to give a general discussion for unique features of mechanics for dynamic behaviors of fractional-order dynamic systems by using the star-coupled network structure as the basic model; and secondly, the general conclusion for the evolution of SMEs' growth is given under the fractional-order world for $\alpha \in (0, 0.8]$, not with true data of SMEs from the market. It is our plan to conduct the research in this direction by applying stochastic models equipped with more complex structures applied to true data of samples on SMEs from the market in the practice.

# Acknowledgements

We like to express our sincere thanks to Professor Hong Ma for his always encouragement and discussion for the study of this work.

## Competing interests

Both authors declare that they have no competing interests.

## Author's contributions

Both authors contributed equally to this work, and read and approved the final manuscript.

## Fundings

This work was supported in part by the National Natural Science Foundation of China with the project numbers: U1811462, 71971031, 12102369; and the National Major Common Key Technology Project of the Ministry of Science and Technology of China with the project number:2018YFB1403005.